\documentclass[preprint,12pt]{elsarticle}

\usepackage{exscale,amssymb}
\usepackage{path}

\usepackage{amsfonts,amssymb,amsmath}

\usepackage{graphicx}
\usepackage{epsfig}

\newtheorem{theorem}{Theorem}

\newtheorem{thm}[theorem]{Theorem}
\newtheorem{lemma}[theorem]{Lemma}

\newcommand{\wt}{\widetilde}

\renewcommand{\leq}{\leqslant}
\renewcommand{\geq}{\geqslant}
\renewcommand{\le}{\leqslant}
\renewcommand{\ge}{\geqslant}
\renewcommand{\phi}{\varphi}

\renewcommand{\theta}{\vartheta}
\DeclareMathOperator{\sign}{sign}

\journal{Linear Algebra Appl.}

\begin{document}

\begin{frontmatter}



\author{Federico Greco\fnref{label2}}
\author{Bruno Iannazzo\fnref{label2}}
\author{Federico Poloni\fnref{label3}}
\fntext[label2]{Dipartimento di Matematica e Informatica, Universit\`a di Perugia, via Vanvitelli, 1, 06123 Perugia, Italy. \texttt{\{greco,bruno.iannazzo\}@dmi.unipg.it}}
\fntext[label3]{Scuola Normale Superiore, piazza dei Cavalieri, 7, 56126 Pisa, Italy. \texttt{f.poloni@sns.it}}

\title{The Pad\'e iterations for the matrix sign function and their reciprocals are optimal}

\begin{abstract}
It is proved that among the rational iterations locally converging with order $s>1$ to the sign function, the Pad\'e iterations and their reciprocals are the unique with the lowest sum of the degrees of numerator and denominator.
\end{abstract}


\end{frontmatter}







\section{Introduction}\label{sec:intro}

The function $s(z)=\sign(z)$ is defined for a nonimaginary complex
number $z$ as the nearest square root of unity. Let $A$ be a
matrix having no purely imaginary eigenvalues. Since $s(z)$ is
analytic at the eigenvalues of $A$, matrix function theory
\cite{highambook} allows one to define $\sign(A)$. The matrix sign
function is less trivial than its scalar counterpart, for instance
it is not locally constant, and it has important applications,
either direct like the solution of algebraic Riccati equations
\cite{kl} and the treatment of certain quantum chromodynamics
models \cite{qcd} or indirect as a basis to compute other
important matrix functions like the matrix square root, the polar
decomposition of a matrix and the geometric mean of two positive
definite matrices \cite{highambook,im}.

A common way to compute the matrix sign function is through
rational iterations of the form $z_{k+1}=\varphi(z_k)$, for some
rational function $\varphi(z)$ having attractive fixed points
at $1$ and $-1$, since any such iteration converges locally to the
function $\sign(z)$. The prototypical example is Newton's method for $z^2-1=0$, but many other iterations have been proposed. Among them, a very popular family is obtained by using the Pad\'e approximants to $f(\xi)=(1-\xi)^{-1/2}$ and the following characterization
\[
    \mathrm{sign}(z)=\frac z{(z^2)^{1/2}}=\frac z{(1-\xi)^{1/2}},
\]
where $\xi=1-z^2$. 
Let the $(m,n)$ Pad\'e approximant to $f(\xi)$ be $P_{m,n}(\xi)/Q_{m,n}(\xi)$, and $m+n\ge 1$.
The iteration
\begin{equation}\label{eq:pf}
    z_{k+1}=\frac{z_k P_{m,n}(1-z_k^2)}{Q_{m,n}(1-z_k^2)}=:\varphi_{2m+1,2n}
\end{equation}
has been proved to be locally convergent to $1$ and $-1$ with order of convergence $m+n+1$ for $m\ge n-1$ \cite{kl91}. The notation $\varphi_{2m+1,2n}$ introduced here highlights the fact that the numerator and denominator of $\varphi_{2m+1,2n}$ have degree $2m+1$ and $2n$, respectively.

We recall that, for integers $m,n\ge 0$, the $(m,n)$ Pad\'e approximant to a function $h(z)$ is a rational function $p(z)/q(z)$, where $p(z)$ and $q(z)$ are polynomials of degree $m$ and $n$, respectively, such that
\[
	h(z)-\frac{p(z)}{q(z)}=O(z^{m+n+1}).
\]
For an introduction to the Pad\'e approximation see the book~\cite{padebook}.

The iterations~\eqref{eq:pf} have been derived by Kenney and Laub \cite{kl91} and are called {\em Pad\'e family of iterations} or just {\em Pad\'e iterations}; they have been considered also in \cite{fs,higham1,highambook,lz09,zh09} for computing matrix functions or invariant subspaces of a matrix.  Observe that the definition of Pad\'e iterations in \cite{kl91} is slightly different from ours, since we exclude the case $m=n=0$, which yields the trivial iteration $z_{k+1}=\varphi_{1,0}(z_k)=z_k$ being not locally convergent to $1$ and $-1$.

Using the identity
\[
\mathrm{sign}(z)=\frac {(z^2)^{1/2}}z=\frac {(1-\xi)^{1/2}}z,
\]
and the  Pad\'e approximants to  $g(\xi)=(1-\xi)^{1/2}$, a different family of
iterations having attractive fixed points at $1$ and $-1$ is obtained. If $p(z)/q(z)$ is the $(m,n)$ Pad\'e approximant to the function $h(z)$ and $h(0)\ne 0$, then $q(z)/p(z)$ is the $(n,m)$ approximant to $1/h(z)$ (see \cite[Theorem 1.5.1]{padebook}), thus the $(m,n)$  Pad\'e approximant to $g(\xi)$ is $Q_{n,m}(\xi)/P_{n,m}(\xi)$. The iteration
\begin{equation}\label{eq:pf2}
    z_{k+1}=\frac{Q_{n,m}(1-z_k^2)}{z_k P_{n,m}(1-z_k^2)}=:\varphi_{2m,2n+1}(z_k)
\end{equation}
is obtained. We call the iterations \eqref{eq:pf2} {\em reciprocal Pad\'e iterations} or {\em reciprocal Pad\'e family}. The possibility to invert the functions defining the Pad\'e iterations is suggested as well by Laub in \cite{laub91} without further discussions.

Many iterations of interest can be retrieved in the Pad\'e family and its reciprocal: $\varphi_{2,1}$, $\varphi_{3,0}$ and $\varphi_{3,2}$ give Newton's method, the Newton-Schulz iteration and Halley's method for $z^2-1=0$, respectively. Among the Pad\'e family \eqref{eq:pf}, the most common iterations are those with constant denominator (i.e., $n=0$ in \eqref{eq:pf}), and the so-called \emph{principal} Pad\'e iterations, namely, those for which the degrees of the numerator and denominator differ by 1 (i.e., $m=n$ or $m=n-1$ in \eqref{eq:pf}). Similarly, the reciprocal Pad\'e family contains iterations with constant numerator ($m=0$ in \eqref{eq:pf2}) or for which the degrees of the numerator and denominator differ by 1 ($m=n$ or $m=n+1$ in \eqref{eq:pf2}).

Prior to Kenney and Laub, the Pad\'e iterations for $m=n-1$ and the reciprocal Pad\'e iterations for $m=n$ have been derived in a different way by Howland \cite{howland83}; Iannazzo
\cite{iannazzo08} proved that the same iterations obtained by
Howland can be retrieved in the family of root-finding algorithms
sometimes called K\"onig family \cite{buff} or basic family
\cite{kalantari} and attributed \cite{buff} to a paper written
by Schr\"oder in 1870 \cite{sch,schen}. 

The Pad\'e family and its reciprocal family are just two of the
infinite families of rational iterations having the square roots
of unity as attractive fixed points. A rationale for their use
can be given by their interesting properties (see \cite{highambook,kl} for the case of the Pad\'e family; some
analogous properties hold for the reciprocal Pad\'e family). We
show that they also have an optimality property: among all
rational iterations having order of local convergence $s>1$ at $1$
and $-1$, they are the unique iterations such that the sum of the degrees of
the numerator and denominator is minimal. This is a highly
desirable property in terms of computational efficiency: in the
{\em generic} case it is cheaper (in terms of the number of
arithmetic operations required) to evaluate $a(z)/b(z)$ than $\wt
a(z)/\wt b(z)$ when $\deg(\wt a(z))+\deg(\wt b(z))>\deg(a(z))+\deg(b(z))$ and Horner's scheme is applied.

Let $s,m,n$ be nonnegative integers such that $s>1$ and $m+n=2s-1$. Observe that, letting $m,n$ vary, the family $z_{k+1}=\varphi_{mn}(z_k)$ is the union of the Pad\'e family and its reciprocal family; the parities of $m$ and $n$ distinguish one from the other.
The following property of $\varphi_{mn}(z)$  is the main result of the paper and will be proved in the next section. 
\begin{thm}\label{thm:gip}
Let $s>1$ be a nonnegative integer.
The functions $\varphi_{mn}(z)$, for $m=0,1,\ldots,2s-1$ and $n=2s-m-1$  define the unique rational iterations of the kind $z_{k+1}=\varphi(z_k)$ such that
\begin{enumerate}
 \item[\emph{\textbf{O1}}] the iteration converges locally to $1$ and $-1$ with order at least $s$;
 \item[\textbf{\emph{O2}}] for every iteration $w_{k+1}=\wt\varphi(w_k)=\wt a(w_k)/\wt b(w_k)$, with $\wt a(z),\wt b(z)$ polynomials, having order at least $s$ in both $1$ and $-1$, it holds that 
$\deg(\wt a(z))+\deg(\wt b(z))\geq \deg(a(z))+\deg(b(z))$, where $a(z)$ and $b(z)$ are coprime polynomials such that $\varphi(z)=a(z)/b(z)$.
\end{enumerate}
Moreover, the iterations have order exactly $s$ in both $1$ and $-1$.
\end{thm}

We consider just the case $s>1$ for two reasons: first, in matrix functions computation, algorithms based on rational iterations are competitive if they converge fast, that is if they are of order at least $2$; second, if $s=1$ a direct computation shows that the unique iterations satisfying \textbf{O1} and \textbf{O2} are the same as the ones obtained for $s=2$.

It is worth noting that a rational iteration satisfying
\textbf{O1} and \textbf{O2} for some $s>1$ is not necessarily the
one whose iteration function can be evaluated with the minimal
cost. In principle, there can be a special rational function which
does not satisfy \textbf{O2} and can be evaluated with fewer
arithmetic operations. In addition, the same iteration can be evalutated with many different schemes yielding different computational costs, relevant in the matrix case (see \cite[Chapter 4]{highambook}). For the principal Pad\'e iterations, a partial fraction expansion \cite{howland83,kl94} and a continued fraction expansion \cite{kbs94} are known and can be used to devise efficient evaluation schemes as in \cite[Algorithms 4.9-4.10]{highambook}.

For the sake of clarity, we recall some basic definitions regarding iterations of the kind $z_{k+1}=\varphi(z_k)$, where $\varphi(z)$ is a rational function and $z_*$ is a fixed point of $\varphi(z)$, that is, $\varphi(z_*)=z_*$. We say that $z_*$ is an \emph{attractive fixed point} if $|\varphi'(z_*)|<1$; in that case the iteration is locally convergent to $z_*$, that is, any initial value $z_0$ sufficiently close to $z_*$ yields a sequence converging to $z_*$. We say that the iteration {\em converges locally to $z_*$ with order $s>1$} if
there exist $M_1,M_2>0$ such that for $z$ sufficiently close to $z_*$ it holds that $M_1|z-z_*|^s\leq |\varphi(z)-z_*|\leq M_2|z-z_*|^s$. Since $\varphi(z)$ is infinitely many differentiable at $z_*$, this is
equivalent to requiring that
\begin{equation}\label{convequiv}
   \varphi(z_*)=z_*,\quad \varphi'(z_*)=\varphi''(z_*)=\cdots=\varphi^{(s-1)}(z_*)=0, \quad \varphi^{(s)}(z_*)\neq 0.
\end{equation}
In particular, an iteration having order $s>1$ at $z_*$ is locally convergent to $z_*$. Further discussion on this topic can be found in any numerical analysis textbook, for instance \cite{sb}.

\section{Proof of the Theorem \ref{thm:gip}}
The proof of Theorem \ref{thm:gip} is given by some Lemmas.
We first prove that imposing \eqref{convequiv} for a rational function $\varphi(z)=a(z)/b(z)$ is equivalent to imposing some conditions on the polynomials $a(z)$ and $b(z)$ and their derivatives. Then we prove that they can only be satisfied if $\deg(a(z))+\deg(b(z))\geq 2s-1$, with equality only for a unique family of polynomials. Finally, we prove that these unique solutions correspond to the Pad\'e family and its reciprocal family.

\begin{lemma}\label{grecolemma}
Let $s>1$ be an integer, $z_{k+1}=\varphi(z_k)=a(z_k)/b(z_k)$ a rational iteration, and $z_*$ be one of its fixed points (in particular, $b(z_*)\neq 0$). The iteration converges locally to $z_*$ with order at least $s$ if and only if
\begin{equation}\label{grecolemma1}
a^{(k)}(z_*)=z_* b^{(k)}(z_*), \quad\text{for
$k=0,1,\dots,s-1$}.
\end{equation}
If moreover $a^{(s)}(z_*)\neq z_* b^{(s)}(z_*)$ then the order is exactly $s$.
\end{lemma}
\noindent {\em Proof.}
The rational iteration converges locally to $z_*$  with order at least $s$ if and only if for $z$ sufficiently close to $z_*$
\begin{equation}\label{cosa1}
\varphi(z)-z_*=O((z-z_*)^s).
\end{equation}
Since $b(z)$ is bounded in a neighborhood of $z_*$, we may multiply the left-hand side of \eqref{cosa1} by $b(z)$ without changing its convergence behavior, thus obtaining
\begin{equation}\label{cosa2}
a(z)-z_*b(z)=O((z-z_*)^s),
\end{equation}
which in turn is equivalent to \eqref{grecolemma1}. As $b(z_*)\neq 0$, $b(z)$ is bounded away from 0 in a neighborhood of $z_*$, thus we may divide the left-hand side of \eqref{cosa2} by $b(z)$ to reverse the previous step and get \eqref{cosa1}.

If $a^{(s)}(z_*)\neq z_* b^{(s)}(z_*)$, then \eqref{cosa2} does not hold anymore if we replace $s$ with $s+1$, and neither does \eqref{cosa1}, i.e., the convergence order is exactly $s$.
\hfill $\square$ \bigskip

Using conditions \eqref{grecolemma1}, we may prove the following bound on the degrees of $a(z)$ and $b(z)$.
\begin{lemma}\label{polonilemma}
Let $s$ be a positive integer, and $a(z),b(z)$ two polynomials, not both null, such that
\begin{equation}\label{grecoeqs}
\left\{%
\begin{array}{ll}
    a^{(k)}(1)\: =\: b^{(k)}(1), \\
    a^{(k)}(-1)=-b^{(k)}(-1), \\
\end{array}%
\right.\quad k=0,1,\ldots,s-1.
\end{equation}
Then $\deg(a(z)) + \deg(b(z))\geq 2s-1$, with the convention that the degree of the zero polynomial is $-1$.

Moreover, for each pair $(m,n)$ of integers such that $m+n=2s-1$,  $\min(m,n)\geq -1$, there are two polynomials
$a_{mn}(z),b_{mn}(z)$ such that $\deg(a_{mn}(z))=m$, $\deg(b_{mn}(z))=n$
and the conditions \eqref{grecoeqs} hold. The polynomials $a_{mn}(z)$ and $b_{mn}(z)$ are
unique up to a multiplicative factor, and
\begin{equation}\label{condizionibuffe}
 a^{(s)}_{mn}(1)\neq b^{(s)}_{mn}(1),\qquad
 a^{(s)}_{mn}(-1)\neq -b^{(s)}_{mn}(-1).
\end{equation}
\end{lemma}
\noindent {\em Proof.}
First notice that we may impose without loss of generality that $\deg(a(z))\geq \deg(b(z))$ throughout the proof.

We shall prove the result by induction. For $s=1$, the result is
clear. The conditions that we must meet are $a(1)=b(1)$,
$a(-1)=-b(-1)$. Thus $a(z)$ and $b(z)$ cannot be both constant,
and the only possibility with $\deg(a(z))=1$, $\deg(b(z))=0$ is
choosing $b(z)=\gamma$, $a(z)=\gamma z$ for some constant $\gamma\neq 0$, while the only possibility with $\deg(a(z))=2$, $\deg(b(z))=-1$ is choosing $b(z)=0$, $a(z)=\gamma (z-1)(z+1)$ for some constant $\gamma\ne 0$.

Let us suppose that the lemma holds true for a given $\bar s-1$, and prove it for $s=\bar s$. Let us take two polynomials $a(z)$, $b(z)$ such that $a^{(k)}(\pm 1)=\pm b^{(k)}(\pm 1)$ for all $k=1,2,\dots, \bar s-1$. If $b(z)\equiv 0$, then $a(z)$ must be a multiple of both $(z-1)^{\bar s}$ and $(z+1)^{\bar s}$, thus it has degree at least $2\bar s$ and the result holds. If $b(z) \neq 0$, we may apply the inductive hypothesis to their derivatives $a^{(1)}(z)$ and $b^{(1)}(z)$, and obtain $\deg(a^{(1)}(z))+\deg(b^{(1)}(z))\geq 2\bar s-3$. This clearly implies $\deg(a(z)) + \deg(b(z))\geq 2\bar s-1$, since the derivative of a polynomial $p(z)$ has degree $\deg(p(z))-1$ (notice that this relation holds also for constant polynomials $p(z)\equiv c \neq 0$, with our choice $\deg(0)=-1$).

Let us turn now to the equality case; we shall prove the uniqueness first, and the existence thereafter, of the two families of polynomials attaining the minimal degrees.

Let $m,n$ be such that $m+n=2\bar s-1$, and $a_{mn}(z)$ and
$b_{mn}(z)$ be two polynomials with $\deg(a_{mn}(z))=m$,
$\deg(b_{mn}(z))=n$ satisfying \eqref{grecoeqs}. If $n=-1$, then
$a(z)$ must be a polynomial multiple of both $(z-1)^{\bar s}$ and
$(z+1)^{\bar s}$ of degree $2\bar s$, and thus
\begin{equation}\label{costrstupida}
 \begin{aligned}
  a_{2s,-1}(z)&=k(z-1)^{\bar s}(z+1)^{\bar s},\\
  b_{2s,-1}(z)&=0
 \end{aligned}
\end{equation}
for some $k\neq 0$. If $n \neq -1$, then  $a^{(1)}_{mn}(z)$ and $b^{(1)}_{mn}(z)$ satisfy the equality conditions of the lemma with $s=\bar s-1$; thus by the uniqueness result it must be the case that
\begin{equation}\label{eq:daint}
\begin{aligned}
a^{(1)}_{mn}(z)&=k a_{m-1,n-1}(z), \\
b^{(1)}_{mn}(z)&=k b_{m-1,n-1}(z),
\end{aligned}
\end{equation}
for some $k \neq 0$. From \eqref{eq:daint} we get
$a_{mn}(z)=a_{mn}(\pm 1)+k\int_{\pm 1}^z a_{m-1,n-1}(t)dt$ and
$b_{mn}(z)=b_{mn}(\pm 1)+k\int_{\pm 1}^z b_{m-1,n-1}(t)dt$.
Imposing $a_{mn}(\pm 1)=\pm b_{mn}(\pm 1)$, a simple manipulation of the resulting
system gives
\begin{equation}\label{costrpol}
\begin{aligned}
a_{mn}(z)&=k A(z)+\frac 12 k(B(1)-B(-1)-A(1)-A(-1)),\\
b_{mn}(z)&=k B(z)+\frac 12 k(A(1)-A(-1)-B(1)-B(-1)),
\end{aligned}
\end{equation}
with $A(z)$ (resp. $B(z)$) a primitive of $a_{m-1,n-1}(z)$ (resp. $b_{m-1,n-1}(z)$). It is now apparent that the two polynomials are uniquely determined up to the multiplicative constant $k$. From the inductive hypothesis $a^{(\bar s-1)}_{m-1,n-1}(\pm 1)\neq \pm b^{(\bar s-1)}_{m-1,n-1}(\pm 1)$, it follows that $a^{(\bar s)}_{mn}(\pm 1)\neq \pm b^{(\bar s)}_{mn}(\pm 1)$.

On the other hand, one can easily check that the polynomials defined by the formulas  \eqref{costrstupida} and \eqref{costrpol} have degree $\deg(a_{mn}(z))=m$, $\deg(b_{mn}(z))=n$ and satisfy \eqref{grecoeqs}. So said polynomials exist.
\hfill $\square$ \bigskip

\noindent {\em Proof of Theorem~\ref{thm:gip}.}
We prove the theorem for a fixed $s>1$.

We first show that for each $m,n\ge 0$ such that $m+n=2s-1$ the iteration $z_{k+1}=\varphi_{mn}(z_k)$ is the unique rational iteration of the kind $z_{k+1}=a(z_k)/b(z_k)$ such that $\deg(a(z))=m$ and $\deg(b(z))=n$, satisfying \textbf{O2} and whose order of local convergence is exactly $s$ (thus it satisfies \textbf{O1} as well). Then, we show that any rational iteration satisfying \textbf{O1} and \textbf{O2} is of the type $z_{k+1}=\varphi_{mn}(z_k)$ for $m+n=2s-1$.

Let $m$ and $n$ be such
that $m+n=2s-1$, with $m,n\ge0$; then, by Lemma~\ref{polonilemma},
there are two polynomials $a_{mn}(z)$ and $b_{mn}(z)$ such that
$\deg(a_{mn}(z))=m$, $\deg(b_{mn}(z))=n$ satisfying
\eqref{grecoeqs} and \eqref{condizionibuffe}.
Since $s>1$, Lemma~\ref{grecolemma} implies that for
$\psi_{mn}(z):=a_{mn}(z)/b_{mn}(z)$ the iteration
$z_{k+1}=\psi_{mn}(z_k)$ converges locally to $1$ and $-1$ with order exactly $s$, thus $\psi_{mn}$ satisfiies
\textbf{O1}.

On the other hand, consider an iteration function $\psi(z)=
a(z)/b(z)$ providing a sequence converging with order at least $s$. By Lemma~\ref{grecolemma}, it follows that conditions \eqref{grecoeqs} hold
and thus $\deg(a(z))+\deg(b(z)) \geq 2s-1=\deg(a_{mn}(z))+
\deg(b_{mn}(z))$ by Lemma~\ref{polonilemma}. Therefore
the iteration $z_{k+1}=\psi_{mn}(z_k)$ satisfies \textbf{O2}. By
the same lemma, equality holds if and only if $a(z)$ and $b(z)$
differ from $a_{mn}(z)$ and $b_{mn}(z)$ by the same multiplicative factor, i.e., when $\psi(z)$ and $\psi_{mn}(z)$ coincide. Thus, this
is the unique iteration satisfying both \textbf{O1} and \textbf{O2} of the kind $z_{k+1}=a(z_k)/b(z_k)$ with $\deg(a(z))=m$ and $\deg(b(z))=n$.

Now we show that $\psi_{mn}(z)$ coincides with $\varphi_{mn}(z)$. Since the numerator of $\varphi_{mn}(z)$ has degree $m$ and its denominator has degree $n$, it is enough to prove that $z_{k+1}=\varphi_{mn}(z_k)$ satisfies \textbf{O1} (thus \textbf{O2}, in view of Lemma \ref{polonilemma}).

a) Odd $m$, $m=2m_1+1$, $n=2n_1$. Let
$h_{\mu\ell}(\zeta):=P_{\mu\ell}(\zeta)/Q_{\mu\ell}(\zeta)$ be the
$(\mu,\ell)$ Pad\'e approximant to $(1-\zeta)^{-1/2}$. Then, in a neighborhood of $\zeta=0$,
\[
    h_{\mu\ell}(\zeta)-(1-\zeta)^{-1/2}=O(\zeta^{\mu+\ell+1}).
\]
Since $\phi_{mn}(z)=zh_{m_1,n_1}(1-z^2)$ and $m_1+n_1=s-1$, we get
$\phi_{mn}(z)-z(z^2)^{-1/2}=O((1-z^2)^{s})$, for $z$ sufficiently close to $1$ or $-1$.
Since $(1-z^2)^s=O((z-1)^s)$  and
$(1-z^2)^s=O((z+1)^s)$ in a neighborhood of $1$ and $-1$ respectively, it holds that
\begin{equation}\label{eq:dentroilteo}
 \begin{aligned}
  \phi_{mn}(z)-1&=O((z-1)^s), & \textrm{for }z\textrm{ in a neighborhood of } 1,\\
  \phi_{mn}(z)+1&=O((z+1)^s), & \textrm{for }z\textrm{ in a neighborhood of } -1,
 \end{aligned}
\end{equation}
and then $z_{k+1}=\varphi_{mn}(z)$ verifies \textbf{O1}.

b) Even $m$, $m=2m_1$, $n=2n_1+1$. By a reasoning similar to case
a, we use the $(\mu,\ell)$ Pad\'e approximant to
$(1-\zeta)^{1/2}$, say $\wt h_{\mu\ell}(\zeta)$, and
$\phi_{mn}(z)=\wt h_{m_1,n_1}(1-z^2)/z$, thus $\varphi_{mn}(z)$ verifies \eqref{eq:dentroilteo} as well and $z_{k+1}=\varphi_{mn}(z)$ verifies \textbf{O1}.

Thus, for each $m$, $n$ such that $m+n=2s-1$, we have $\varphi_{mn}(z)=\psi_{mn}(z)$ and the iteration $z_{k+1}=\varphi_{mn}(z_k)$ is the unique iteration satisfying \textbf{O1} and \textbf{O2} and whose numerator has degree $m$ and denominator has degree $n$.

Let $w_{k+1}=a(w_k)/b(w_k)$ be a rational iteration satisfying \textbf{O1} and \textbf{O2}. Let $m'=\deg(a(z))$ and $n'=\deg(b(z))$. By Lemmas~\ref{grecolemma} and \ref{polonilemma} one has $m'+n'\ge 2s-1$. We claim that $m'+n'=2s-1$; if on the contrary $m'+n'>2s-1$ then there exist $m\le m'$ and $n\le n'$ such that $m+n=2s-1$ and $z_{k+1}=\varphi_{mn}(z_k)$ satisfies \textbf{O1}, thus $w_{k+1}=a(w_k)/b(w_k)$ cannot satisfy \textbf{O2} and we get a contradiction. Finally, by the aforementioned uniqueness result we conclude that $a(z)/b(z)$ must coincide with $\varphi_{m',n'}(z)$ up to a multiplicative factor.
\hfill $\square$ \bigskip

\section*{Acknowledgement}

We wish to thank Nicholas J. Higham who pointed out that properties
\textbf{O1} and \textbf{O2} do not necessarily imply the minimal evaluating cost and two anonymous referees whose suggestions improved the paper. This work was supported by the grant PRIN 20083KLJEZ of the Italian Ministero dell'Istruzione, dell'Universit\`a e della Ricerca, and by the Istituto Nazionale di Alta Matematica.

\end{document}